\documentclass[12pt]{amsart}

\newtheorem{theorem}{Theorem}[section]
\newtheorem{proposition}[theorem]{Proposition}
\newtheorem{lemma}[theorem]{Lemma}
\newtheorem{corollary}[theorem]{Corollary}

\theoremstyle{definition}
\newtheorem{definition}[theorem]{Definition}
\newtheorem{example}[theorem]{Example}

\theoremstyle{remark}
\newtheorem{remark}[theorem]{Remark}

\numberwithin{equation}{section}

\usepackage{color}		% Need the color package
\usepackage{epsfig}
\usepackage{pinlabel} 
\begin{document}

\title{$K$-stable equivalence for knots in Heegaard surfaces} 
\author{Alice Stevens}
\address{Department of Mathematics, University of California, Davis, California 95616}
\email{astevens@math.ucdavis.edu}
\thanks{Research partially supported by NSF grants.}
%\subjclass{57N10}
\date{February 20, 2009.}
\begin{abstract}
Let $K$ be a knot embedded in a Heegaard surface $S$ for a closed orientable 3-manifold $M.$ We define $K$-stable equivalence between pairs $(S, K)$ and $(S', K)$ in $M$, and we prove that any two pairs are $K$-stably equivalent in $M$ if they have the same surface slope.
\end{abstract}

\maketitle

%%%%%%%%%%%%%%%%%%%%%%%%%%%%%%%%%%%%%%%%%%%%%%%%%%%
%                                           INTRODUCTION                                                                                                          %  
%%%%%%%%%%%%%%%%%%%%%%%%%%%%%%%%%%%%%%%%%%%%%%%%%%%

\section{Introduction}
Embeddings of knots in interesting surfaces of 3-manifolds are relevant to both 3-manifold theorists and knot theorists. An interesting surface that exists in any compact orientable 3-manifold $M$ is a \emph{Heegaard surface}, which decomposes $M$ into two simple homeomorphic pieces $V$ and $W$, called \textit{handlebodies}. We denote this decomposition $M=V\cup_{S} W$. \emph{Torus knots}, which are knots that can be embedded in genus one Heegaard surfaces for $S^3$, are well understood. \emph{Double torus knots}, which can be embedded in genus two Heegaard surfaces for $S^3$, were studied in \cite{Hi} by Hill, and in \cite{HM} by Hill and Murasugi. Morimoto \cite{Mo} studied the \textit{h-genus} of a knot in $S^3$, which is the minimal genus of any Heegaard surface for $S^3$ in which the knot can be embedded. Knots in Heegaard surfaces also appear in Dehn surgery theory. In \cite{Be}, Berge studied a special subfamily of double torus knots, called \emph{doubly primitive knots}, or \emph{Berge knots}, which admit lens space surgeries. In \cite{De}, Dean studied \textit{twisted torus knots}, which admit small Seifert fibered Dehn surgeries. This paper is part of a project to study questions of equivalence between knots in Heegaard surfaces for a closed orientable 3-manifold $M$.

%%%%%%%%%%%%%%%%%%%%%%%%%%%%%%%%%%%%%%%%%%%%%%%%%%%

Let $S$ be a Heegaard splitting surface for a compact orientable 3-manifold $M$, and let $K$ be a knot embedded in $S$. We call $(S, K)$ a \emph{$K$-splitting pair for $M$}, and we call the Heegaard splitting $M=V \cup_{(S, K)}W$ a \textit{$K$-splitting for M}. Intuitively speaking, two splitting pairs $(S, K)$ and $(S^{'}, K^{'})$ are \textit{equivalent} in $M$ if we can push one pair onto the other in $M$. More formally, we have the following definition:

\begin{definition} \label{e}
(Equivalent splitting pairs) Let $M$ be a compact orientable 3-manifold. If $(S, K)$ and $(S^{'}, K^{'})$ are splitting pairs for $M$, then $(S, K)$ and $(S^{'}, K^{'})$ are $\emph{equivalent}$ in $M$ if there is an ambient isotopy of $M$ that maps $(S, K)$ onto $(S^{'}, K^{'})$.
\end{definition}

When are two splitting pairs equivalent, and how many different equivalence classes exist? Definition \ref{e} implies four obvious conditions that must be satisfied if $(S, K)$ and $(S^{'}, K^{'})$ are equivalent in $M$. We must have that $S$ and $S'$ are isotopic as Heegaard surfaces of $M$, $K$ and $K'$ are of the same knot type, $S-K$ and $S^{'}-K^{'}$ have the same number of connected components, and the surface slope of $K$ with respect to $S$ must be equal to the surface slope of  $K^{'}$ with respect to $S'$.

\begin{figure}[htb]
\centering
\includegraphics[width=2.7in]{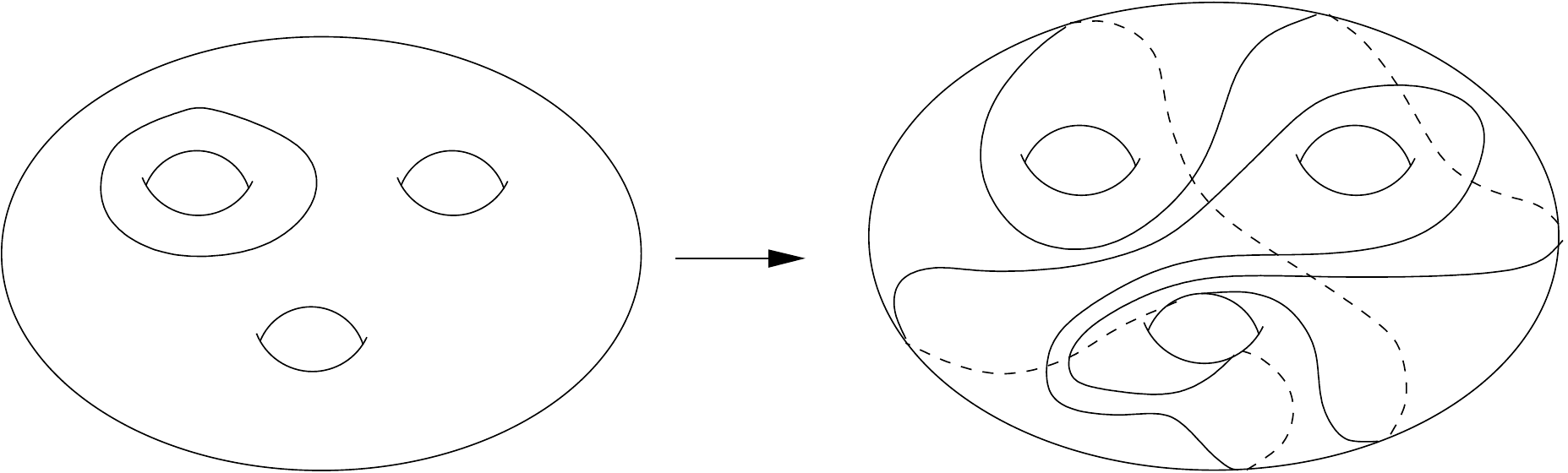}
\caption{Equivalent $K$-splitting pairs in $S^3$}\label{fig:two}
\end{figure}

Even in the case that $M$ is $S^3$ and $K$ is the unknot, it may be difficult to determine whether two $K$-splitting pairs are equivalent in $M$. Consider, for example, the \textit{equivalent} embeddings of the unknot in Figure~\ref{fig:two}. In both splitting pairs, the unknot is embedded in a genus three Heegaard surface for $S^3$ as a non-separating curve with surface slope zero. 

\begin{figure}[htb]
\centering
\includegraphics[width=1in]{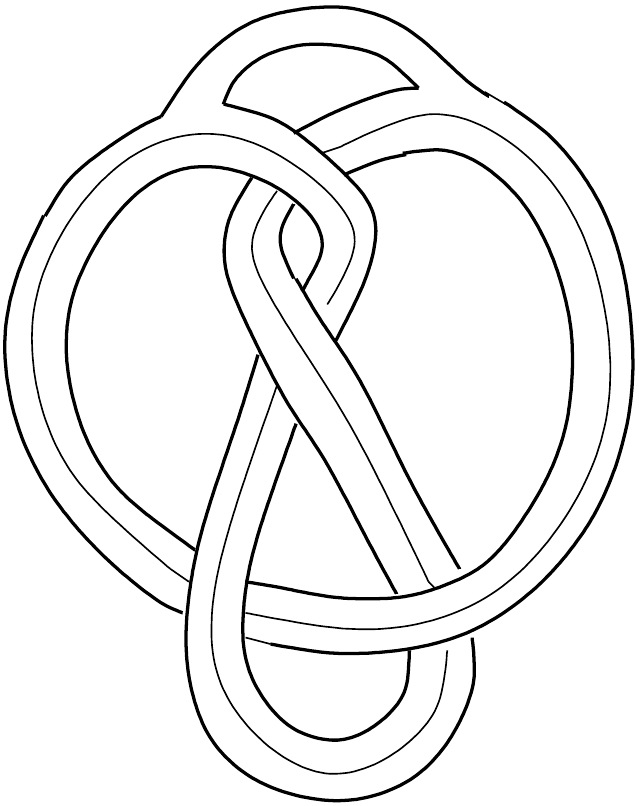}
\caption{The figure eight knot embedded in a genus two Heegaard surface for $S^3$}\label{eight}
\end{figure}

There \textit{are} inequivalent $K$-splitting pairs in $S^3$ that satisfy the four conditions mentioned above. For example, let $K$ be a tunnel number one knot. There is a natural way to construct an embedding of $K$ as a non-separating curve in a genus two Heegaard surface for $S^3$ so that it has any given surface slope (see Figure ~\ref{eight}). Note that $K$ can be pushed slightly into the interior of the handlebody so that it is a core of a handle, and the meridian disk corresponding to the tunnel is uniquely defined. If $K$ has two non-isotopic tunnels, we can construct two $K$-splitting pairs for $K$. If these $K$-splitting pairs are equivalent in $S^3$, then there is an isotopy of $S^3$ that maps one pair onto the other. By pushing $K$ slightly into the handlebody, we obtain an isotopy of $S^3-K$ that maps one tunnel onto the other, a contradiction. 

%%%%%%%%%%%%%%%%%%%%%%%%%%%%%%%%%%%%%%%%%%%%%%%%%%%

A \emph{stabilization} of a Heegaard splitting surface is the process of adding an unknotted tube to the surface. This can be characterized more formally as follows: a Heegaard splitting $V\cup_{S}W$ is \emph{stabilized} if and only if we can find $D_1\subset V$ and $D_2 \subset W$  such that $\lvert D_1 \cap D_2 \rvert= 1.$ Two Heegaard surfaces for a compact orientable 3-manifold $M$ are \textit{stably equivalent} if they have a common stabilization. The Reidemeister-Singer theorem states that any two Heegaard surfaces for a compact orientable 3-manifold are \textit{stably equivalent}. 

\begin{figure}[htb]
\centering
\includegraphics[width=2.5in]{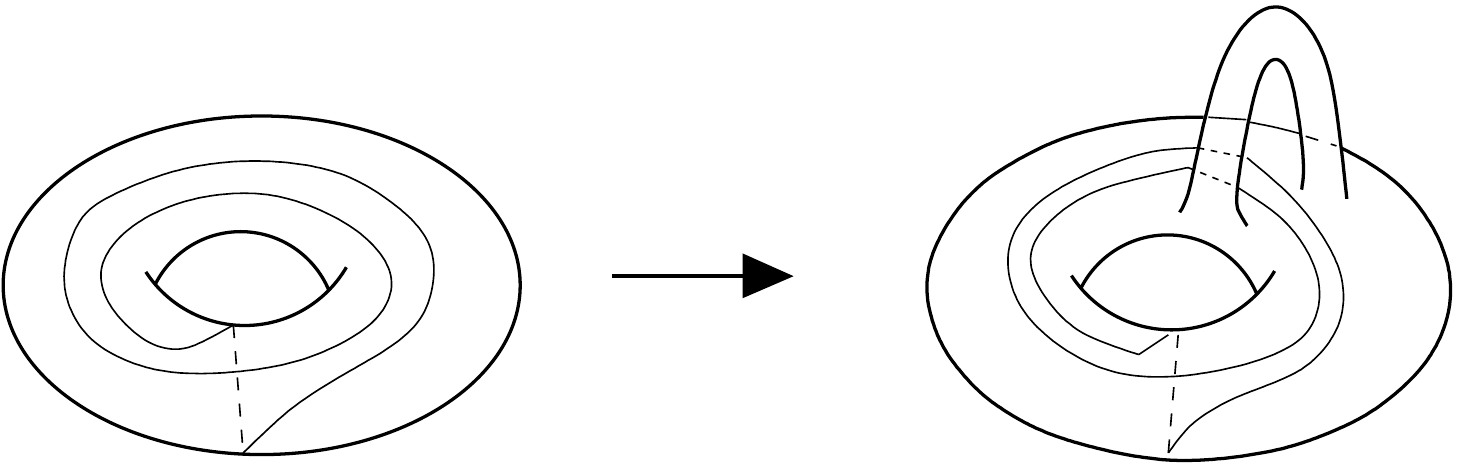}
\caption{$K$-stabilization}\label{fig:three}
\end{figure}

A similar notion of equivalence can be defined for $K$-splitting pairs. Roughly speaking, we define a \textit{$K$-stabilization} of the $K$-splitting pair $(S, K)$ to be the addition of an unknotted tube to the surface $S-K$, where the tube may straddle the knot as in Figure~\ref{fig:three}. The formal definition, which follows, is analogous to the definition of stabilization found in \cite{Sc}:

\begin{definition} ($K$-stabilization) \label{ks}
Suppose $V \cup_{(S, K)}W$ is a $K$-splitting for a closed orientable 3-manifold $M$. Let $\alpha$ be a properly embedded arc in $W$ parallel to an arc $\beta$ in $S$, and such that $\partial\alpha \cap K=\emptyset$. Add a neighborhood of $\alpha$ in $W-K$ to $V$, and remove it from $W$. This adds a 1-handle to each handlebody, creating two handlebodies $\tilde{V}$ and $\tilde{W}$ of genus one greater than $V$ and $W$. We say that $\tilde{V} \cup_{\tilde{S}} \tilde{W}$ is a \textit{$K$-stabilization} of $V \cup_{(S, K)} W.$ 
\end{definition}

$K$-stabilization can also be characterized as follows: a $K$-splitting $M=V \cup_{(S, K)} W$ is \textit{$K$-stabilized} if and only if there are properly embedded disks $D_1\subset V$, and $D_2 \subset W$ such that $\vert \partial D_1 \cap \partial D_2 \vert = 1$ and $\vert D_1\cap K \vert = 0$ (the disk $D_2$ can intersect $K$ multiple times). Note that the standard notion of stabilization is a special case of $K$-stabilization. We now have the following notion of equivalence for $K$-splitting pairs:

\begin{definition} ($K$-stable equivalence)
If $(S, K)$ and $(S^{'}, K)$ are two $K$-splitting pairs for a compact orientable 3-manifold $M$, then $(S, K)$ and $(S^{'}, K)$ are $\emph{$K$-stably equivalent}$ in $M$ if they have a common $K$-stabilization.
\end{definition}

The goal of this paper is to show that two $K$-splitting pairs for $M$ are \textit{$K$-stably equivalent} as long as they have the same surface slope. Our main theorem is the following:

\begin{theorem} \label{main} 
Let $K$ be a knot in a closed orientable 3-manifold $M.$ Suppose $(S, K)$ and $(S^{'}, K)$ are two $K$-splittings for $M$ such that $K$ is embedded in both surfaces with surface slope $m$, then $(S, K)$ and $(S^{'}, K)$ are $K$-stably equivalent.
\end{theorem}

%%%%%%%%%%%%%%%%%%%%%%%%%%%%%%%%%%%%%%%%%%%%%%%%%%%

In order to prove the theorem, one might attempt to apply the Reidemeister-Singer theorem directly; however, the presence of the knot is an obstruction to this approach. The theorem is proved using \emph{weak reduction} and \emph{amalgamation} adapted to the study of knots in Heegaard surfaces. Given any $K$-splitting pair for $M$, we can $K$-stabilize to enable a decomposition of $M$ into three $(K$-)splittings: a $K$-splitting of a solid torus $\mathcal{T}$ that is isotopic to $\eta(K)$, a $K$-splitting of the product manifold $T^2 \times [0, 1]$, and a Heegaard splitting of the knot complement $M-\mathcal{T}$. This decomposition allows us to isolate the knot and apply the classical Reidemeister-Singer Theorem to the knot complement. As a result, it is possible to construct an isotopy of $M$ that maps a $K$-stabilization of $(S,K)$ onto a $K$-stabilization of $(S^{'}, K)$. 

The outline of this paper is as follows: Section \ref{a} contains preliminary definitions and propositions, including a discussion of \textit{surface slope} in Subsection \ref{b}, and a description of \textit{$K$-weak reduction} and \textit{$K$-amalgamation} in Subsection \ref{tools}. Section \ref{lemmas} contains a sequence of four lemmas used in the proof of Theorem \ref{main}. Section \ref{proof} contains the proof of Theorem \ref{main}. 

I would like to thank my thesis advisor, Jennifer Schultens, for her patience, support, and encouragement, as well as for many helpful suggestions and conversations. I would also like to thank Abigail Thompson and Kei Nakamura. The work in this paper is part of my thesis; the motivation is to understand equivalence of embeddings of knots in Heegaard surfaces for 3-manifolds.

%%%%%%%%%%%%%%%%%%%%%%%%%%%%%%%%%%%%%%%%%%%%%%%%%%%
%                                              PRELIMINARIES                                                                                                       %  
%%%%%%%%%%%%%%%%%%%%%%%%%%%%%%%%%%%%%%%%%%%%%%%%%%%

\section{Preliminaries}\label{a}

\subsection{Basic definitions}

\noindent For standard definitions and facts about 3-manifolds, see \cite{He} and \cite{Ja}. For standard definitions and facts about knots, see \cite{Ro}.

\begin{definition} \label{cb} (Compression body)
Let $F$ be a closed, orientable, possibly disconnected surface, and let $\mathcal{O}$ be a collection of 3-balls. Construct a 3-manifold $V$ by attaching 1-handles to the disjoint union of $F \times [0, 1]$ and $\mathcal{O}$, along $F \times \{1\} \subset F \times [0, 1]$ and $\partial \mathcal{O}$, in such a way that the resulting manifold is connected and orientable. Any 3-manifold homeomorphic to one constructed in this manner is called a $\emph{compression body}$. The boundary component $F \times \{0\}$ is denoted $\partial_{-}V$, and $\partial V-\partial_{-}V$ is denoted $\partial_{+}V.$ When $V$ is a handlebody, $\partial_{-}V$ is defined to be empty. 
\end{definition}

\begin{figure}[htb]
\centering
\includegraphics[width=3in]{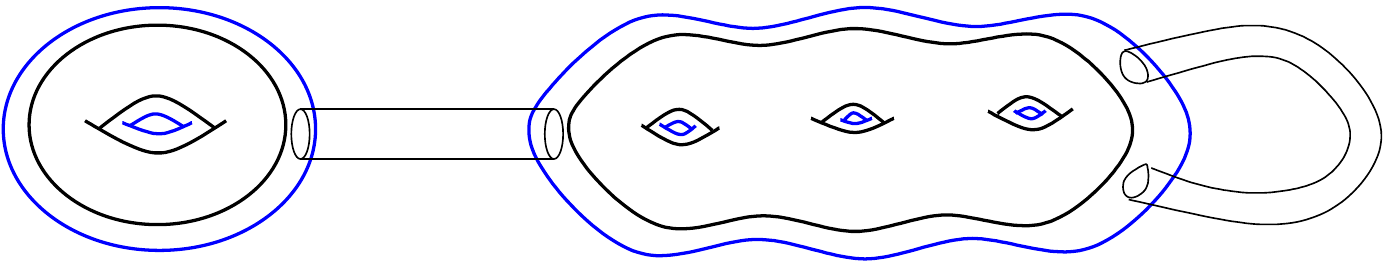}
\caption{A compression body} \label{fig:box27}
\end{figure}

For the compression body in Figure~\ref{fig:box27}, $F$ is the disjoint union of a torus and a genus three surface. We have added two 1-handles to $F\times\{1\}$.
   
Let $M$ be a \textit{compact} orientable 3-manifold. A \textit{Heegaard splitting of M} is a decomposition $M=V \cup_{S} W,$ where $V$ and $W$ are compression bodies, and $S$ is a closed, connected, orientable, embedded surface satisfying $S=V \cap W=\partial_{+}V=\partial_{+}W.$ $S$ is called a $\emph{Heegaard splitting surface}.$ If $K$ is a knot embedded in $S$, then $(S, K)$ is a $K$-$splitting$ $pair$ for $M$.

\vspace{.05in}

Let $K$ be a knot embedded in a 3-manifold $M.$ We denote a regular neighborhood of $K$ in $M$ by $\eta(K)$. In general, $\eta(\cdot)$ will be used to denote a regular neighborhood. A $\emph{tunnel system}$ for a knot $K$ in $S^3$ is a collection of disjoint arcs $\tau=\{\tau_i\}$ such that $K\cap\tau=\partial \tau$ and $S^3-int(\eta(K\cup\tau))$ is a handlebody. If $n$ is the minimal number of arcs in any tunnel system for $K$, then $K$ is said to have \emph{tunnel number n}, denoted $t(K)=n.$

%%%%%%%%%%%%%%%%%%%%%%%%%%%%%%%%%%%%%%%%%%%%%%%%%

\subsection{Surface slope} \label{b}

Let $M$ be a 3-manifold, and let $S$ be a surface in the interior of $M$. If $K$ is embedded in $S$, then $\partial \eta(K)$ is a torus, and $\partial \eta(K) \cap S$ consists of two curves, $\alpha_1$ and $\alpha_2$. The isotopy class of these curves in $\partial \eta(K)$ is called the \textit{surface slope of $K$ with respect to $S$}.

If $K$ is a knot in a Heegaard surface for $M$, we may simply refer to $m$ as the \textit{surface slope of the $K$-splitting pair $(S, K)$}. In the canonical basis on $\partial\eta(K),$ the surface slope can be identified by a fraction. Note that the surface slope is always integral, since the $\alpha_i$ are isotopic to a core of $\partial \eta(K)$. 

\begin{example}
 Let $K$ be a separating curve in a Heegaard surface $S$ for $S^3$, and suppose $(\lambda, \mu)$ is the canonical basis for $\partial \eta(K)$, i.e. $\lambda$ is homologically trivial in $S^3 - int(\eta(K))$ and $\mu$ generates $H_1(S^3 - int(\eta(K))).$ Then the surface slope of $K$ with respect to $S$ is 0, since the two components of $S-K$ are Seifert surfaces for $K$. 
\end{example}

The following lemma allows us to compute the surface slope for nonseparating curves in surfaces in $S^3.$

\begin{lemma} \label{compute}
Let $(S, K)$ be a $K$-splitting pair for $S^3.$ Suppose $(\mu, \lambda)$ is the canonical basis for $\partial \eta (K)$. If $\alpha_1 \sqcup \alpha_2$ is the two component link in $S^3$ given by $\partial{\eta(K)} \cap S$, then the surface slope of $K$ with respect to $S$ is equal to the linking number $lk(\alpha_1, \alpha_2)$.
 \end{lemma}
 
\begin{proof}
Let $m$ be the surface slope of $K$ with respect to $S.$ Then with respect to the canonical basis, $[\alpha_i]=[\lambda] + m[\mu]$, and so $\lambda$ and $\alpha_i$ differ by $m$ meridianal twists. As a preferred longitude, $\lambda$ satisfies $lk(K, \lambda)=0,$ and therefore $lk(K, \alpha_1)=m$. But $K$ is the core of $\eta(K)$, and $\alpha_2$ can be isotoped to $K$ in $\eta(K)$, so $lk(\alpha_1, \alpha_2)=lk(\alpha_1, K)=m$. 
\end{proof}

The surface slope of a splitting pair $(S, K)$ is invariant under isotopy of the ambient manifold, i.e., if two splitting pairs are equivalent, then they must have the same surface slope. This fact allows us to create infinitely many inequivalent $K$-splitting pairs, which is the subject of the Proposition \ref{inf} below. 

Recall that the $h$-genus of a knot $K$, denoted $h(K),$ is the smallest genus of any Heegaard surface for $S^3$ in which the knot can be embedded. 

\begin{proposition} \label{inf}
Let $K$ be a knot in $S^3$. Then for any $m \in \mathbb{Z}$, there is a Heegaard surface of genus less than or equal to $h(K)+1$ in which $K$ can be embedded as a non-separating curve with surface slope $m$. 
\end{proposition}

\begin{proof}
Let $m \in \mathbb{Z}$. We will perform a connect sum of triples. Let $V_1 \cup_{(S_1, K)} W_1$ be the $K$-splitting of $S^3$ which realizes the $h$-genus of $K$, and suppose the surface slope of $K$ with respect to $S_1$ is $n$. Let $B_1$ be a 3-ball in $S^3$ such that $B_1 \cap K$ is a single trivial arc in $B_1$, and $B_1 \cap S_1$ is a disk, so $B_1$ intersects $V_1$ and $W_1$ in a single ball.

Let $K'$ denote the unknot, and let $V_2 \cup_{(S_2, K')} W_2$ be the genus one $K'$-splitting of $S^3$ such that $K'$ is embedded in $S_2$ as a curve that wraps once in the longitudinal direction, and $m-n$ times in the meridianal direction. By Lemma \ref{compute}, the surface slope of $V_2 \cup_{(S_2, K')} W_2$ is $m-n.$ Let $B_2$ be a 3-ball in $S^3$ such that $B_2 \cap K$ is a single trivial arc in $B_2$, and $B_2$ intersects $S_2$ in a disk. 
 
Finally, let $V=cl(V_1-B_1)\cup cl(V_2-B_2)$ and $W=(W_1-B_1)\cup cl(W_2-B_2)$ where $\partial B_1 \cap V_1$ is identified with $\partial B_2 \cap V_2,$ and $\partial B_1 \cap W_1$ is identified with $\partial B_2 \cap W_2$, and $\partial(B_1 \cap K)$ is identified with $\partial(B_2 \cap K')$. Then $V \cup_{(S_1 \# S_2, K)} W$ is a genus $h(K)+1$ $K$-splitting of $S^3$ with $K$\#$K'=K$ embedded in $S_1$\#$S_2$ having surface slope $n+(m-n)=m$. 
\end{proof}

\begin{proposition} \label{merid}
Let $K$ be a knot in $S^3.$ Suppose $V \cup_{(S, K)}W$ is a genus $g$ $K$-splitting for $S^3$. If there is a non-separating meridian disk $D$ embedded in $V$ or $W$ such that $D$ intersects $K$ in one point, then for any $m \in \mathbb{Z}$ there is a genus $g$ $K$-splitting pair with surface slope $m$.
\end{proposition}

\begin{figure}[htb]
\centering
\includegraphics[width=1.1in]{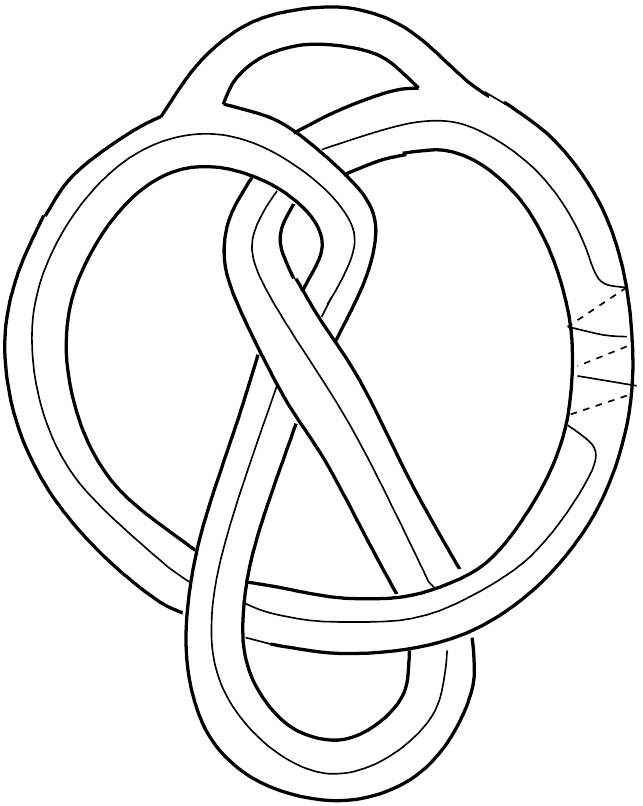}
\caption{Dehn twist on a non-separating meridian disk} \label{fig:box40}
\end{figure}

\begin{proof} 
By Lemma \ref{compute}, each Dehn twist performed on the surface $S$ along $\partial D$ modifies the surface slope of $(S, K)$ in $M$ by $\pm 1$, without changing the knot type of $K$. Compare Figure~\ref{eight} with Figure~\ref{fig:box40}. 
\end{proof}

The \textit{tunnel number} of a knot $K$ in $S^3$ gives bounds for the $h$-genus: if $t(K)$ is the \textit{tunnel number} of $K$, then $t(K) \leq h(K) \leq t(K) + 1$ (see \cite{Mo}). Next, we apply Proposition \ref{merid} to the case that $h(K)=t(K)+1$.

\begin{corollary}
If $h(K)=t(K)+1$, then $K$ has infinitely many distinct minimal $h$-genus $K$-splitting pairs, in particular, for any $m \in \mathbb{Z}$, there exists a Heegaard surface of genus $h(K)$ in which $K$ is embedded as a non-separating curve with surface slope $m.$
\end{corollary}

\begin{proof} 
Let $\Gamma=K \cup \tau$, where $\tau$ is a tunnel system for $K$. Then $K$ is a core of a handle of the genus $t(K)+1$ handlebody $\eta(\Gamma)$, so we can isotope $K$ into $\partial \eta (\Gamma)$. Since $K$ was a core of a handle of $\Gamma,$ a meridian disk of this handle intersects $K$ in one point. Now we apply Proposition \ref{merid} to modify this embedding to have any desired surface slope. See Figure~\ref{fig:box40}.
\end{proof}

If $K$ is a knot embedded in an arbitrary compact orientable 3-manifold $M$, then the above constructions apply. As a result, for any $m \in \mathbb{Z}$ there is a $K$-splitting pair for $M$ with surface slope $m$.

%%%%%%%%%%%%%%%%%%%%%%%%%%%%%%%%%%%%%%%%%%%%%%%%%%%

\subsection{$K$-stable equivalence}

Recall that a \textit{$K$-stabilization} is the process of adding an unknotted tube to a $K$-splitting $V \cup_{(S, K)} W$, and that if
two $K$-splitting pairs $(S, K)$ and $(S^{'}, K^{'})$ have a $common$ $K$-$stabilization$, we say that they are \textit{$K$-stably equivalent}.

In the following lemma, we prove an equivalent characterization for $K$-stabilization, mentioned in the introduction of this paper. The proof of the lemma follows the proof of Lemma 3.1 in \cite{Sc}.

\begin{lemma}
The $K$-splitting $M=V \cup_{(S, K)} W$ is \textit{$K$-stabilized} if and only if there are properly embedded disks $D_1\subset V$, and $D_2 \subset W$ such that $\vert \partial D_1 \cap \partial D_2 \vert = 1$ and $\vert D_1\cap K \vert = 0$.  
\end{lemma}

\begin{proof}
If $(S, K)$ is $K$-stabilized, then according to Definition \ref{ks} we can find a disk $D$ that is bounded by $\alpha \cup \beta$. If we let $D_2=D$, and let $D_1$ be a cocore of $\eta(\alpha),$ we have found two disks satisfying the above conditions. 
Conversely, one can compress $S$ along $D_1$ to obtain $S_1$, which bounds a handlebody. Since $\partial D_1 \cap K = \emptyset$, $K$ is embedded in $S_1$. Compressing $S$ along $D_2$ yields a surface $S_2$ which also bounds a handlebody (although $K$ will not be embedded in $S_2$). 

Next we would like to show that $S_1$ bounds a handlebody on \textit{both} sides, and therefore it is a Heegaard surface for $M$. A neighborhood of $D_1$ and $D_2$, $\eta(D_1 \cup D_2)$, is a ball. The boundary of this ball intersects $S_1$ in one hemisphere $H^{+}$, and intersects $S_2$ in the other hemisphere $H^{-}$. This allows us to isotope the surface $S_2$ to $S_1$. Since $S_2$ bounded a handlebody, we see that $S_1$ bounds a handlebody on \emph{both} sides. Therefore, $S_1$ is a Heegaard splitting surface for $M$, and since $K$ is embedded in $S_1$, $(S_1, K)$ is a $K$-splitting pair for $S^3.$ If $\alpha$ is taken to be the core of a 1-handle dual to $D_1,$  then $\partial\alpha \cap K=\emptyset$, and we may $K$-stabilize to obtain our original $K$-splitting pair $(S, K)$.
\end{proof}

\begin{figure}[htb]
\centering
\includegraphics[width=2.5in]{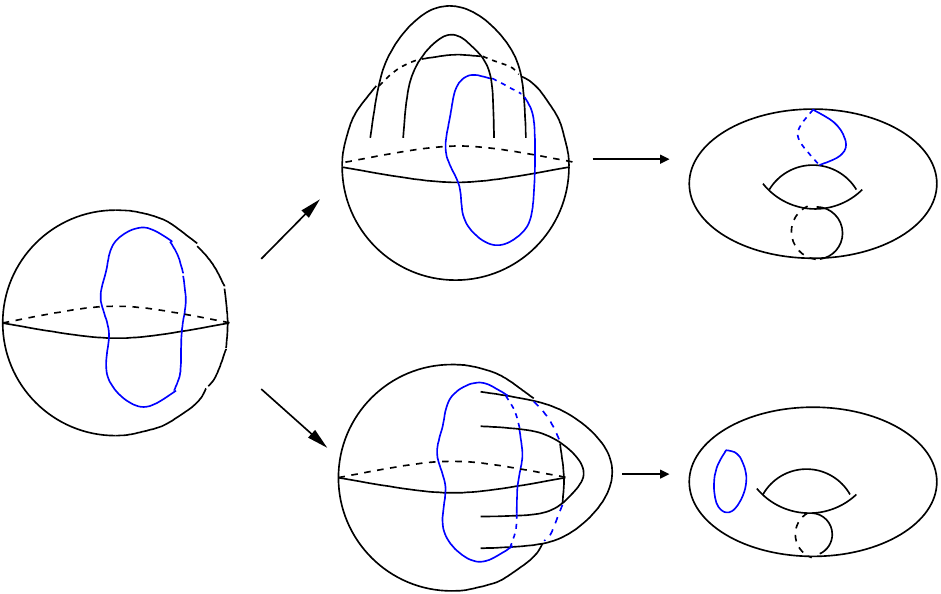}
\caption{$K$-stabilization is not unique}\label{fig:notthesame.pdf}
\end{figure}

Unlike stabilization, $K$-stabilization is not unique. Figure \ref{fig:notthesame.pdf} shows two inequivalent ways to $K$-stabilize the $K$-splitting pair $(S, K)$, where $K$ is the unknot and $S$ is the 2-sphere in $S^3$.

An important fact about $K$-stabilization and surface slope is the following:

\begin{proposition}
The surface slope of a knot is invariant under $K$-stabilization.
\end{proposition}

\begin{proof}
Suppose $M=V \cup_{(S, K)} W$ is a $K$-splitting. Recall that when we $K$-stabilize, we remove a  tubular neighborhood of an arc $\alpha$, where $\alpha$ is properly embedded in one of the handlebodies, say $W$, and $\eta(\alpha)\cap\eta (K) =\emptyset $. The resulting surface $\tilde{S}$, is just $\partial (V \cup \eta(\alpha))$, and therefore $S \cap \partial n(K)=\tilde{S}\cap \partial n(K),$ so the two $K$-splitting pairs $(S, K)$ and $(\tilde{S}, K)$ have the same surface slope. This is sufficient to show that the slope is invariant under $K$-stabilization, since the surface slope is invariant under ambient isotopy: any isotopy of $M$ restricts to an isotopy of the submanifold $\partial \eta(K)$. The curves $S \cap \partial \eta (K)$ will be sent to the same isotopy class in $\partial \eta(K)$. In conclusion, if the surface slopes of two $K$-splitting pairs are different, then there is no hope for a common $K$-stabilization.
\end{proof}

%%%%%%%%%%%%%%%%%%%%%%%%%%%%%%%%%%%%%%%%%%%%%%%%%%%
%                                           K-Weak reduction and K-amalgamation                                                                 %  
%%%%%%%%%%%%%%%%%%%%%%%%%%%%%%%%%%%%%%%%%%%%%%%%%%%

\subsection{$K$-weak reduction and $K$-amalgamation} \label{tools}

A $K$-splitting can be decomposed into a collection of Heegaard splittings such that $two$ of the splittings are $K$-splittings. The following construction follows from the standard theory of $weakly$ $reducible$ Heegaard splittings and their corresponding $induced$ Heegaard splittings. See \cite{CG} and \cite{Sch1}.

\begin{definition} ($K$-weakly reducible) A $K$-splitting $M=V \cup_{(S, K)} W$ of a compact orientable 3-manifold $M$ will be called \textit{$K$-weakly reducible} if there are two essential properly embedded disks $D_1 \subset V$, and $D_2 \subset W$, such that $\partial D_1 \cap \partial D_2 =\emptyset$ and neither disk intersects the knot $K.$ We also consider \emph{collections} of $K$-weak reduction disks $\Delta_V \subset V$ and $\Delta_W \subset W$, where for any $D_i \subset \Delta_V$ and any $D_j \subset \Delta_W,$ $\partial D_i \cap \partial D_j = \emptyset$, and each disk is disjoint from $K$.
\end{definition}

If there is a collection of $K$-weak reduction disks $\Delta_V \cup \Delta_W$ for the $K$-splitting pair $(S, K)$ such that both $\Delta_V$ and $\Delta_W$ are non-empty, then one can compress $S$ simultaneously into both handlebodies along $\Delta_V \cup \Delta_W$, leaving $K$ embedded in exactly one of the components of the resulting surface $S^{*}$. We would like to describe the connected components of $M-S^{*}$.

After compressing $S$ along $\Delta_V$, the handlebody $V$ is split into $\overline{V}=V-\eta(\Delta_V)$. We will denote the connected components of $\overline{V}$ by $\overline{V_i}$. Compressing along $\Delta_{W}$ is equivalent to attaching 2-handles to the components $\overline{V_i}$. We will denote the result of attaching the 2-handles to $\overline{V_i}$ by $C_i$.

Symmetrically, compressing $S$ along $\Delta_{W}$ splits the handlebody $W$ into $\overline{W}=W-\eta(\Delta_W)$, which is a collection of connected components $\overline{W_i},$ to which we attach the 2-handles $\eta(\Delta_V)$. The set of all components $C_i$ for $i=1,...n,$ is the set of connected components of $M-S^{*}$. 

Each component $C_i$ has a Heegaard splitting induced by $V\cup_{(S, K)}W$. There are different ways to describe this induced Heegaard splitting, but we will follow the construction in the proof of Lemma 2.4 of \cite{Sch1}: choose one component $C_i$, and without loss of generality, assume that $C_i=\overline{V_i} \cup \eta(\Delta_{W^{'}})$,  where $\Delta_{W^{'}} \subset \Delta_W$, and $\eta(\Delta_{W^{'}})$ is the collection of 2-handles attached to $\overline{V_i}$. Consider the fattened up version of $C_i$, $C_i^*=C_i \cup (\partial C_i \times [0, 1])$. We now construct a Heegaard splitting for $C_i^{*}$. Let $V_i=\overline{V_i}$, and $W_i=( \partial C_i \times [0, 1]) \cup (1$-$handles)$, where the 1-handles are dual to the 2-handles $\eta(\Delta_{W^{'}}$). Then $V_i \cup W_i$ is a Heegaard splitting for $C_i^{*}$. Note that the Heegaard splitting of $C_i^{*}$ is a Heegaard splitting of $C_i$ as well, and that the collection of boundary components of $\cup_{i=1}^{n}\partial C_i^{*}$ is $S^{*}$.

Now suppose $S_i^{*}$ is one of the connected components of the compressed surface $S^{*}.$ If $K \subset S_i^{*},$ then $K$ and $S_i^{*}$ are incident to two of the components of $M-S^{*}$, say $C_i$ and $C_{j}.$ By construction, there is a copy of the knot $K$ in both of the induced Heegaard surfaces for $C_i^{*}$ and $C_j^{*}$. There is also a copy of $K$ in $\partial C_i^{*}$ and in $\partial C_j^{*}$. The two copies of $K$ in $C_i^{*}$ cobound an annulus in the compression body $(\partial C_i \times [0, 1]) \cup (1$-$handles)$, and the same is true in the case of $C_j^{*}$. This annulus records information about the way the knot is embedded in the Heegaard surface, and this information will be useful when we attempt to reconstruct our original $K$-splitting pair. 
We have now decomposed our original $K$-splitting into the family of splittings $V_1 \cup_{S_1},...,V_i \cup_{(S_i, K)}W_i,...,V_j \cup_{(S_j, K)}W_j,...,V_n\cup_{S_n}W_n$. There is a way to recover the original splitting, and this will be discussed next.
\vspace{.08in}

%%%%%%%%%%%%%%%%%%%%%%%%%%%%%%%%%%%%%%%%%%%%%%%%%%%

$Amalgamation$ was first formalized by Schultens. A rigorous treatment can be found in \cite{Sch1}. For the remainder of this Subsection, we assume $M$ and $M_i$ are compact orientable 3-manifolds. 

The intuitive idea of amalgamation is as follows: let $M_1$ and $M_2$ be 3-manifolds, and let $R$ be a connected surface such that $R\subset \partial M_i.$ $M_1$ and $M_2$ can be glued together along $R$ by a homeomorphism to create a 3-manifold $M.$ Given any pair of Heegaard splittings $M_1=V_1 \cup_{S_1} W_1$ and $M_2=V_2 \cup_{S_2}W_2$, one can construct a Heegaard splitting for $M$.

\begin{figure}[ht!]
\labellist
\small \hair 2pt
\pinlabel $V_1$ at 176 330
\pinlabel $W_1$ at 176 249
\pinlabel $V_2$ at 176 162
\pinlabel $W_2$ at 176 76
\pinlabel $S_1$ at 311 288
\pinlabel $R$ at 274 205
\pinlabel $S_2$ at 311 122
\endlabellist
\centering
\includegraphics[width=1.5in]{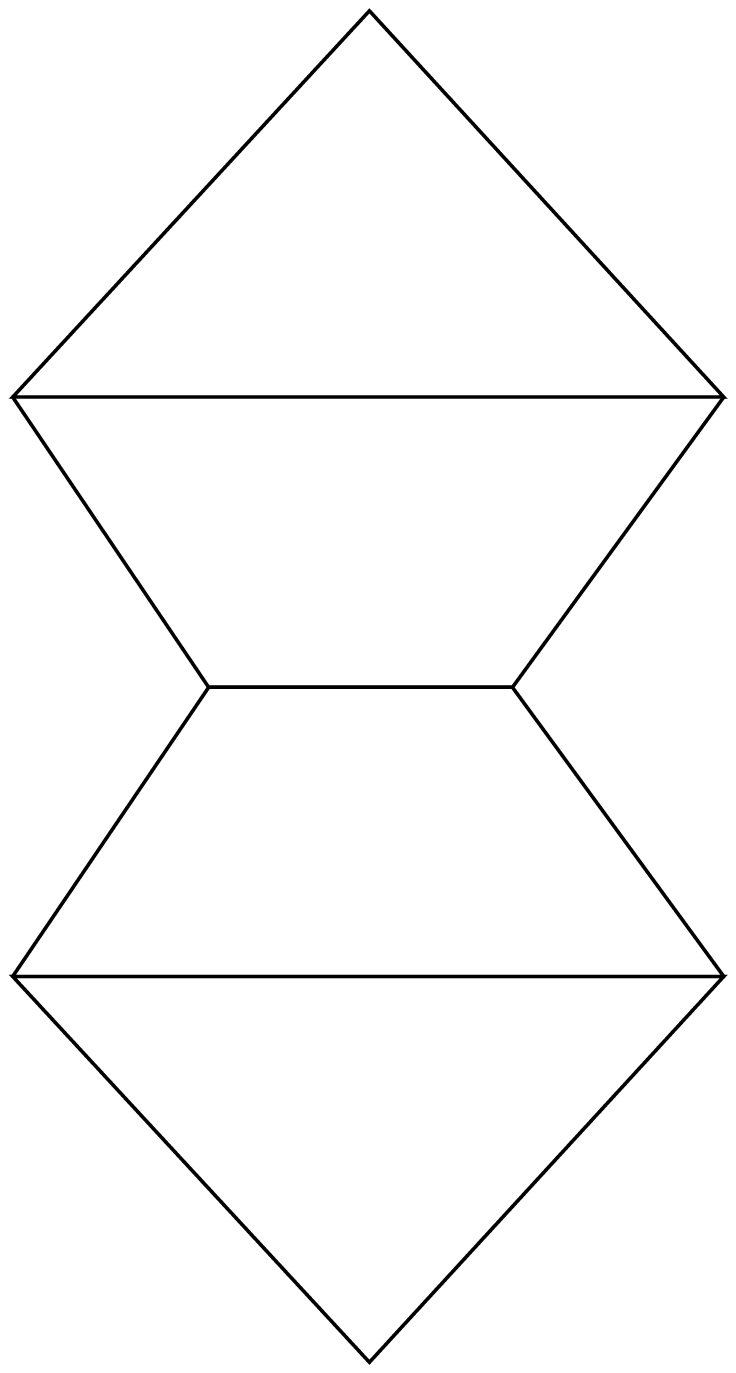}
\caption{Schematic before amalgamation is performed} \label{fig:box12}
\label{fig:am}
\end{figure}

First, we assume that $R$ is one of the boundary components of the compression bodies $W_1$ and $V_2$ (see Figure~\ref{fig:box12} for a schematic where $R=\partial_{-}W_1=\partial_{-}V_2$). Given a handle structure for $W_1$, let $Q_1$ denote the collection of 1-handles that are attached to the copy of $R\times I$ in $W_1$, and given a handle structure for $V_2$, let $Q_2$ denote the set of 1-handles of $V_2$ that are attached to the copy of $R\times I$ in $V_2$.

\begin{figure}[ht!]
\labellist
\small \hair 2pt
\pinlabel $Q_1$ at 300 280
\pinlabel $Q_2$ at 225 50
\pinlabel $S_1$ at -33 240
\pinlabel $R$ at -37 175
\pinlabel $S_2$ at -33 100
\pinlabel $\rightarrow$ at 605 168
\pinlabel $S$ at 1218 175
\endlabellist
\centering
\includegraphics[width=3.5in]{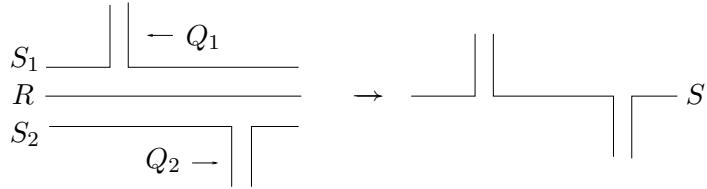}
\caption{A schematic for amalgamation} \label{fig:box13}
\end{figure}

By a small isotopy, we can guarantee that the attaching disks $\mathcal{D}_1$ of $Q_1$ are disjoint from the attaching disks $\mathcal{D}_2$ of $Q_2$ in $R$. Next, collapse both product structures $R \times [0, 1] \subset W_1, V_2$ into the surface $R$ (see Figure~\ref{fig:box13}). What remains is a Heegaard splitting $S$ of $M$ composed of the two compression bodies $V_1 \cup Q_2$ and $W_2 \cup Q_1$.

The following proposition states that the Heegaard splitting that results from amalgamation is independent of any choices made during the construction, such as the choice of handle structure for the compression bodies. 

\begin{proposition} \label{amalgunique}
The operation of amalgamation is well defined.
\end{proposition}

\begin{proof}
See \cite{La} or \cite{Sch2}.
\end{proof}

Amalgamation can also be defined for a $K$-splitting $V_1 \cup_{(S_1, K)} W_1$ and a standard Heegaard splitting $V_2 \cup_{S_2} W_2$ along some surface $R \subset \partial_{-}W_1, \partial_{-}V_2$. However, unlike amalgamation of Heegaard surfaces, this kind of amalgamation is not uniquely defined. Choices of handle structure or isotopy can lead to inequivalent amalgamations.

Under certain conditions, amalgamation involving a $K$-splitting pair \textit{is} unique. For example, recall from Subsection \ref{tools} that after performing a $K$-weak reduction we obtain $K$-splittings for the two components $C_i^{*}$ and $C_j^{*}$. By construction, copies of $K$ are embedded in $\partial C_i^{*}$ and in $\partial C_j^{*}$, as well as in the Heegaard surfaces of these components (we can find two annuli that are cobounded by copies of $K$). If we amalgamate these $K$-splittings along $\partial C_i^{*}$ and $\partial C_j^{*}$ by identifying the copies of the knot $K$, then we obtain a unique $K$-splitting pair. We call this \textit{$K$-amalgamation}. 
																			      
%%%%%%%%%%%%%%%%%%%%%%%%%%%%%%%%%%%%%%%%%%%%%%%%%%%

\begin{proposition} \label{hi}
Suppose the $K$-splitting $M=V \cup_{(S,K)} W$ is the result of amalgamating the $K$-splittings $M_1=V_1 \cup_{(S_1, K)} W_1$ and $M_2=V_2 \cup_{(S_2, K)} W_2$ along a surface $R \subset \partial_{-}W_1, \partial_{-}V_2$. Suppose also that both $K$-splittings contain a copy of $K$ in $R$ that cobounds an annulus with $K \subset S_i$. Then the $K$-amalgamation $M=V \cup_{(S,K)}W$ is unique. 
\end{proposition}

\begin{proof}
If we identify the copies of the knot $K$ in $R$, then the proof is essentially the same as the proof of the main result in \cite{Sch2}.
\end{proof}

%%%%%%%%%%%%%%%%%%%%%%%%%%%%%%%%%%%%%%%%%%%%%%%%%%%

\begin{figure}[ht!]
\labellist
\small \hair 2pt
\pinlabel $(S,K)$ at 60 200
\pinlabel $V$ at 353 269
\pinlabel $W$ at 353 152
\pinlabel $V_1$ at 940 344
\pinlabel $(S_1,K)$ at 1125 290
\pinlabel $W_1$ at 940 256
\pinlabel $\leftarrow K\subset R$ at 1098 207
\pinlabel $V_2$ at 940 161
\pinlabel $$ at 1093 134
\pinlabel $W_2$ at 940 82
\pinlabel $(S_2,K)$ at 1128 114
\endlabellist
\centering
\includegraphics[width=4in]{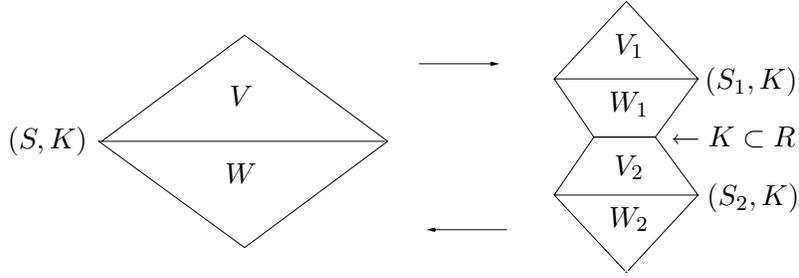}
\caption{$K$-weak reduction and $K$-amalgamation for $n=2$} \label{fig:box34}
\end{figure}

Proposition 2.8 from \cite{Sch1} states that weak reduction and amalgamation are essentially inverse operations. The proof relies on the fact that amalgamation is uniquely defined (Proposition \ref{amalgunique}). Using Proposition \ref{hi}, we can state this result for $K$-splitting pairs:

\begin{lemma} \label{inverses}
Let $\Delta_V \cup \Delta_W$ be a $K$-weak reducing collection of disks for $M=V\cup_{(S, K)}W$, and suppose that $V_1 \cup_{S_1} W_1,...,V_i \cup_{(S_i, K)} W_i,...,V_j \cup_{(S_j, K)} W_j,..., V_n \cup_{S_n} W_n$ are the induced $($K-$)$splittings. Then $V \cup_{(S, K)} W$ is the $(K$-$)$amalgamation of $V_1 \cup_{S_1} W_1,...,V_i\cup_{(S_i, K)} W_i,...,V_j \cup_{(S_j, K)} W_j,...,V_n \cup_{S_n} W_n$. Figure~\ref{fig:box34} is a schematic for n=2. 
\end{lemma}

\begin{proof}
The proof follows from [Sch1, Proposition 2.8] and Proposition \ref{hi} above.
\end{proof}

%%%%%%%%%%%%%%%%%%%%%%%%%%%%%%%%%%%%%%%%%%%%%%%%%%%
%                                           FOUR LEMMAS                                                                                                           %  
%%%%%%%%%%%%%%%%%%%%%%%%%%%%%%%%%%%%%%%%%%%%%%%%%%%

\section{Four lemmas} \label{lemmas}

Not all $K$-splitting pairs $(S, K)$ for a closed orientable 3-manifold $M$ are $K$-weakly reducible, for example, consider a torus knot embedded in a genus one splitting of $S^3$. The first two lemmas of this section specify $K$-stabilizations for a $K$-splitting pair in $M$ that result in a $K$-weakly reducible splitting. The final two lemmas of this section apply the ideas of Subsection \ref{tools} to this $K$-weakly reducible splitting. The four lemmas should be read sequentially, and each subsequent lemma should be thought of as building on the previous lemma. For the remainder of this section, we assume $M$ is a closed orientable 3-manifold. First, we will need the following definition:

\begin{definition}
Suppose $K$ is a knot in $M$ and $\mathcal{T}$ is a solid torus in $M$ that is isotopic to $\eta(K)$, and such that $K$ is a longitude of $\partial \mathcal{T}$. Then we will call $\mathcal{T}$ a $\emph{collar}$ of $K$ in $M$. Note that collars are not unique, and there is one collar per isotopy class of longitudes.
\end{definition}

\noindent We now state the first lemma, which specifies a $K$-stabilization that allows us to `peel' the knot off of the Heegaard surface. 

\begin{lemma} \label{one} 
Let $V \cup_{(S, K)} W$ be a genus $g$ $K$-splitting pair of $M$ where $K$ has surface slope $m$ with respect to $S$. Then there is a $K$-stabilization $\tilde{V} \cup_{(\tilde{S}, K)} \tilde{W}$, and a disk $D_1\subset \tilde{V}$ such that $D_1\cap K=\emptyset$, and $\tilde{V}-D_1$ consists of two components: a solid torus that is a collar of $K$, and a genus $g$ handlebody. The surface slope of $K$ with respect to the boundary of the collar is $m$.
\end{lemma}

\begin{proof} 
Consider the following $K$-stabilization: let $\eta(K) \cong S^1 \times D^2 $ be a neighborhood of $K$ in $M$ and $D={\{p\}} \times D^2$ for some $p \in K$. Let $\alpha=\partial D \cap V$. Then $\alpha$ is a properly embedded arc in $V$, $\partial \alpha \subset S-K,$ and $\alpha$ is parallel to the arc $\beta=D \cap S$ in $S$ as required in the definition of $K$-stabilization. Denote this $K$-stabilized Heegaard splitting by $\tilde{V} \cup_{(\tilde{S}, K)}\tilde{W}$. Note that $K$ is now a non-separating curve in $\partial \tilde{V}=\partial \tilde{W}$, since $\partial(D \cap \tilde{V})$ is a loop intersecting $K$ in one point.

Now let $D_1=\partial \eta(K) \cap \tilde{V}.$ Then $D_1$ is a properly embedded separating disk in $\tilde{V}$. The component of $\tilde{V}-D_1$ containing $K$ is isotopic to the solid torus $\eta(K)$. This solid torus is a collar of $K$ that contains $K$ in its boundary with the same surface slope that $K$ had in $S.$ The other component of $\tilde{V}-D_1$ is a genus $g$ handlebody that is isotopic to $V$.
\end{proof}

The next lemma allows us to find a second disk $D_2$, transforming the original $K$-splitting pair into a $K$-weakly reducible splitting pair.

\begin{lemma} \label{two}
There is a second $K$-stabilization $\hat{V} \cup_{(\hat{S}, K)} \hat{W}$, and a properly embedded disk $D_2$ in $\hat{W}$, such that $D_1\cap D_2 =\emptyset$, $D_2 \cap K =\emptyset$ and $\hat{W}-D_2$ consists of two components: a solid torus that is a collar of $K$, and a genus $g+1$ handlebody. The surface slope of $K$ with respect to the boundary of the collar is again $m$.
\end{lemma}

\begin{proof}
The second $K$-stabilization follows the same procedure as the $K$-stabilization in Lemma \ref{one}, with some slight modifications. Let $\eta^{'}(K) \cong S^1\times D^2$ be a neighborhood of $K$ such that $\eta^{'}(K) \subset \eta(K)$, i.e., $\eta^{'}(K)$ is strictly contained in $\eta(K)$ from Lemma \ref{one}. Let $D^{'}={\{p\}} \times D^2$ for some $p \in K$. Let $\beta=\partial D^{'} \cap W$. Then $\beta$ is a properly embedded arc in $W$, $\partial \beta \subset S-K,$ and $\beta$ is parallel to the arc $D^{'} \cap S$ in $S$ as required. Denote this new Heegaard splitting $\hat{V} \cup_{(\hat{S}, K)} \hat{W}$.

Now let $D_2=\partial \eta^{'}(K) \cap \hat{W}.$ Then $D_2$ is a properly embedded separating disk for $\hat{W}$. The component of $\hat{W}-D_2$ containing $K$ is a solid torus isotopic to $\eta^{'}(K).$ This solid torus is a collar of $K$ that contains $K$ in it's boundary with the same surface slope that $K$ had in $S$. The other component of $\hat{W}-D_2$ is a genus $g+1$ handlebody that is isotopic to $\tilde{W}$.
\end{proof}

%%%%%%%%%%%%%%%%%%%%%%%%%%%%%%%%%%%%%%%%%%%%%

Note that the disks $D_1$ and $D_2$ from the above lemmas are disjoint, and that neither disk intersects the knot, so $\hat{V} \cup_{(\hat{S}, K)} \hat{W}$ is $K$-weakly reducible splitting pair. The next lemma describes the effect of compressing $(\hat{S}, K)$ along $D_1$ and $D_2$.

\begin{lemma} \label{three!}
Compressing $\hat{S}$ along $D_1$ and $D_2$ results in the surface $S^{*}=T_1^2 \sqcup T_2^2 \sqcup\Sigma_g$, where the $T_i^2$ are tori, and $\Sigma_g$ is a genus $g$ surface. $M-S^{*}$ is a decomposition of $M$ into 4 components: $C_1$ is a solid torus, $C_2=T^2 \times [0, 1]$, $C_3$ is a genus $g+1$ handlebody with a 2-handle attached along $\partial D_1$ so that $\partial C_3=T^2 \sqcup \Sigma_g$, and finally, $C_4$ is a genus $g$ handlebody. 
\end{lemma}

\begin{proof}
First, we describe the surface $S^{*}$ that results from compressing $\hat{S}$ along $D_1 \cup D_2$. Since $D_1 \cup D_2$ is a weak $K$-reducing pair of disks, we may compress $\hat{S}$ simultaneously into $\hat{V}$ and $\hat{W}$ to obtain a surface $S^{*}$.

Compressing along the separating disk $D_1$ yields a genus two surface and a genus $g$ surface, $\Sigma_g$. Compressing along the separating disk $D_2$ cuts the genus two surface into two tori $T_1^2$ and $T_2^2,$ so $S^{*}=T_1 \sqcup T_2 \sqcup\Sigma_g$.

We now describe the four components of $M-S^{*}$, noting that each component of $M-S^{*}$ is a component of $\hat{V}-D_1$ possibly with the 2-handle $\eta(D_2)$ attached, or a component of $\hat{W}-D_2$ possibly with the 2-handle $\eta(D_1)$ attached, as described in Section \ref{tools}. One component of $\hat{V}-D_1$ is a genus $g$ handlebody $C_4$, while the other component is a genus 2 handlebody. The 2-handle $\eta(D_2)$ is attached to the genus 2 handlebody along $\partial D_2$ yielding $T^2 \times [0, 1]$ which we denote $C_2$. One component of $\hat{W}-D_2$ is a solid torus we will denote $C_1$, and the other is a genus $g+1$ handlebody. The 2-handle $\eta(D_1)$ is attached to the genus $g+1$ handlebody along $\partial D_1$, creating a compression body $C_3$.
\end{proof}

%%%%%%%%%%%%%%%%%%%%%%%%%%%%%%%%%%%%%%%%%%%%%%%

The final lemma describes the induced Heegaard splittings of the components $C_i^{*}$ from Lemma \ref{three!}.

\begin{lemma} \label{four}
There is a decomposition of $M$ into three components: a solid torus $\mathcal{T}$ (that is collar of $K$), a product manifold $T^2 \times [0, 1]$, and the knot complement $M-\mathcal{T}$, each having a Heegaard splitting induced by $\hat{V} \cup_{(\hat{S}, K)} \hat{W}$. Furthermore, the original $K$-splitting $\hat{V}\cup_{(\hat{S}, K)} \hat{W}$ can be obtained from the $(K$-$)$amalgamation of these three splittings. 
\end{lemma}

\begin{proof}
Recall the components $C_i$, for $i=1,..., 4$ from Lemma \ref{three!}. The Heegaard splitting $M=\hat{V} \cup_{(\hat{S}, K)} \hat{W}$ induces a Heegaard splitting on each $C_i^{*}$, where $C_i^{*}$ is the thickened version of $C_i$ as described in Subsection \ref{tools}. We note that the knot $K$ lies in $\partial C_1$, and $\partial C_2$ (as well as in the Heegaard splitting surfaces of these components).
Following the construction in Subsection \ref{tools}, we now describe the four induced Heegaard splittings $C_i^{*}=V_i \cup_{S_i} W_i$.

\noindent
\begin{enumerate}
\item $C_1^*$: $V_1$ is a solid torus, and $W_1$ is the product manifold $T^2 \times [0, 1]$. This is the trivial Heegaard splitting of a solid torus, so $S_1=\partial_{+} V_1=\partial_{+} W_1$ is a torus, as is $\partial_{-} W_1.$ A copy of $K$ is embedded in the tori $S_1$ and $\partial C_1^*$, as a longitude with surface slope $m$. 

\item $C_2^*: V_2$ is a genus two handlebody, and $W_2=(\partial C_2 \times I) \cup (1$-$handle)$ is a compression body with $S_2=\partial_{+}V_2=\partial_{+} W_2$ a genus two surface, and $\partial_{-} W_2=T^2 \sqcup T^2$. A copy of the knot $K$ is embedded in $S_2$ with slope $m$, as well as in \textit{one} of the tori of $\partial_{-} W_2=T^2 \sqcup T^2 \subset \partial C_2^*$. 

\item  $C_3^*: V_3$ is a genus $g+1$ handlebody, and $W_3=(\partial C_3 \times I) \cup (1$-$handle)$ is a compression body with $S_3=\partial_{+}V=\partial_{+} W_3$ a genus $g+1$ surface, and $\partial_{-} W_3=T^2 \sqcup \Sigma_g$ where $\Sigma_g$ is a genus $g$ surface.

\item $C_4^*$: $V_4$ is a genus $g$ handlebody, and $W_4$ is $\Sigma_g \times [0, 1]$, where $\Sigma_g$ is a genus $g$ surface. This is the trivial Heegaard splitting of a genus $g$ handlebody, so $S_4$ is a genus $g$ surface, and $\partial_{-} W_4$ is also a genus $g$ surface.
\end{enumerate}

\noindent Amalgamate the Heegaard splittings $C_3^{*}$ and $C_4^{*}$ along the genus $g$ surface $\Sigma_g$ in $\partial_{-} W_3=T^2 \sqcup \Sigma_g$, and $\partial_{-} W_4 = \Sigma_g$ (see Figure~\ref{fig:box29}) The result of this amalgamation is a Heegaard splitting $\hat{V_3} \cup_{\hat{S_3}} \hat{W_3}$ of the complement of $K,$ $M-\mathcal{T}$, where $\hat{V_3}$ is a compression body with $\partial_{-} \hat{V_3}=T^2$, $\hat{W_3}$ is a genus $g+1$ handlebody, and finally $\hat{S_3}=S_3$ is a genus $g+1$ surface. 

\begin{figure}[htb]
\centering
\includegraphics[width=4in]{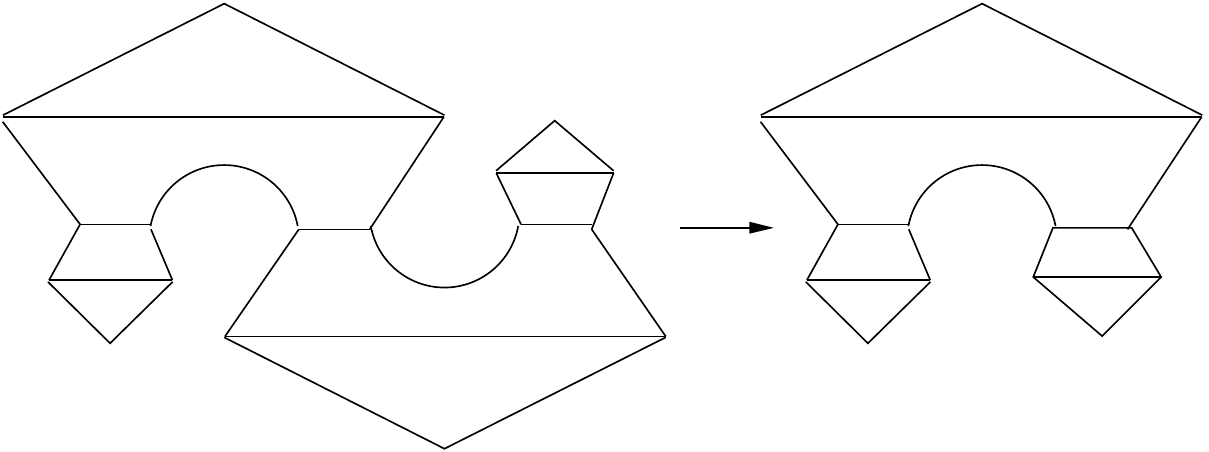}
\caption{Schematic of the three Heegaard splittings} \label{fig:box29}
\end{figure}

By Lemma \ref{inverses}, if we $K$-amalgamate the Heegaard splittings $V_1 \cup_{(S_1, K)} W_1$ and $V_2 \cup_{(S_2, K)} W_2$ along the torus boundary they share, and we amalgamate $V_2 \cup_{(S_2, K)} W_2$ and $\hat{V_3} \cup_{\hat{S_3}} \hat{W_3}$ along the torus boundary they share (there is not copy of $K$ in these tori), the result will be equivalent to the original doubly $K$-stabilized Heegaard splitting $\hat{V} \cup_{(\hat{S}, K)} \hat{W}$, completing the proof of the lemma.
\end{proof}

\begin{remark} \label{rmk1}
Recall that for any $K$-splitting of $M$, the solid torus $\mathcal{T}$ of the previous lemma is isotopic to $\eta(K)$ in $M$. We note that by construction, if the surface slope of $(S, K)$ is $m$, then the $K$-splitting $\mathcal{T}=\hat{V}_1 \cup_{(\hat{S}_1, K)} \hat{W}_1$ obtained in Lemma \ref{four} is unique up to equivalence of splitting pairs in $\mathcal{T}$. The same is true of the $K$-splitting of the product manifold $T^2 \times [0, 1].$ 
\end{remark}

%%%%%%%%%%%%%%%%%%%%%%%%%%%%%%%%%%%%%%%%%%%%%%%%%%%
%                                           PROOF OF MAIN THEM                                                                                              %
%%%%%%%%%%%%%%%%%%%%%%%%%%%%%%%%%%%%%%%%%%%%%%%%%%%

\section{Proof of the Main Theorem} \label{proof}

\noindent The goal of this section is to prove Theorem \ref{main}, which we now restate:
\vspace{.05in}

\noindent $\mathbf{Theorem}$ $\mathbf{\ref{main}}$ 
\emph{Let $K$ be a knot in a closed orientable 3-manifold $M.$ Suppose $(S, K)$ and $(S^{'}, K)$ are two $K$-splittings for $M$ such that $K$ is embedded in both surfaces with surface slope $m$, then $(S, K)$ and $(S^{'}, K)$ are $K$-stably equivalent}.
 
\begin{proof}

Apply Lemmas \ref{one} and \ref{two} to $V \cup_{(S, K)} W$ and to $V^{'} \cup_{(S^{'}, K)} W^{'}$ in order to obtain the $K$-weakly reducible splittings $M=\hat{V}\cup_{(\hat{S}, K)} \hat{W}$ and $M=\hat{V}^{'}\cup_{(\hat{S}^{'}, K)} \hat{W}^{'}.$ By Lemma \ref{four}, we can decompose each of these splittings into three induced splittings $\hat{V_i}\cup_{(\hat{S_i}, K)} \hat{W_i}$ and $\hat{V_i}^{'}\cup_{(\hat{S_i}^{'}, K)} \hat{W_i}^{'}$, for $i=1, 2, 3$. Recall that $i=1$ corresponds to a $K$-splitting for a solid torus $\mathcal{T}$ that is isotopic to $\eta(K)$ in $M$, $i=2$ corresponds to a $K$-splitting for a product manifold $T^2 \times [0, 1],$ and $i=3$ corresponds to a Heegaard splitting of the knot complement $M-\mathcal{T}$. 

By construction, the induced $K$-splittings of the solid tori are equivalent as pairs (see Remark \ref{rmk1}), as are the $K$-splittings for the product manifold $T^2\times [0, 1]$. The induced Heegaard splittings $\hat{S_2}$ and $\hat{S_2}^{'}$ of the knot complements may not be isotopic in $M-\mathcal{T}$, but by the Reidemeister-Singer theorem, they are \textit{stably equivalent} in $M-\mathcal{T}$. We can assume the stabilizations were performed before applying lemmas \ref{one} through \ref{four}, and by abuse of notation we refer to the stabilized surfaces as $\hat{S_2}$ and $\hat{S_2}^{'}$ once again. 

By Lemma \ref{inverses}, the $K$-splitting pairs for $M$ obtained by ($K$-)amal-
gamation of $\hat{V_i}\cup_{\hat{S_i}} \hat{W_i}$ for $i= 1, 2, 3$, and by ($K$-)amalgamation of $\hat{V_i}^{'}\cup_{\hat{S_i}^{'}} \hat{W_i}^{'}$ for $i=1, 2, 3$ are equivalent as $K$-splitting pairs. Thus, the original pairs ${(S, K)}$ and $(S^{'}, K)$ are $K$-stably equivalent. 
\end{proof}

\bibliographystyle{amsalpha}

\end{document}